\documentclass{article}
\usepackage[utf8]{inputenc}
\usepackage[english]{babel}
\usepackage{amsthm}
\usepackage{amsmath}
\usepackage{graphicx}
\usepackage{url}
\usepackage{hyperref}
\usepackage{amsfonts}

\begin{document}

\title{On the complex magnitude of Dirichlet beta function}

\author{Artur Kawalec}

\date{}
\maketitle

\begin{abstract}
In this article, we derive an expression for the complex magnitude of the Dirichlet beta function $\beta(s)$ represented as a Euler prime product and compare with similar results for the Riemann zeta function. We also obtain formulas for $\beta(s)$ valid for an even and odd $k$th positive integer argument and present a set of generated formulas for $\beta(k)$ up to $11$th order, including Catalan's constant and compute these formulas numerically. Additionally, we derive a second expression for the complex magnitude of $\beta(s)$ valid in the critical strip from which we obtain a formula for the Euler-Mascheroni constant expressed in terms of zeros of the Dirichlet beta function on the critical line. Finally, we investigate the asymptotic behavior of the Euler prime product on the critical line.
\end{abstract}

\section{Introduction}
The Riemann zeta function is classically defined by a series

\begin{equation}\label{eq:1}
\zeta(s) = \lim_{k\to\infty}\sum_{n=1}^{k}\frac{1}{n^s},
\end{equation}
which is absolutely convergent for $\Re(s)>1$, where $s$ is a complex variable $s=\sigma+it$. For the first few positive integers, the values of the Riemann zeta function are

\begin{equation}\label{eq:9}
\begin{aligned}
\zeta(1) &\sim\gamma+\log(k) \\
\zeta(2) &=\frac{\pi^2}{6} \\
\zeta(3) &=1.20205690315959\dots \\
\zeta(4) &=\frac{\pi^4}{90} \\
\zeta(5) &=1.03692775514337\dots.
\end{aligned}
\end{equation}
For $s=1$, the series diverges asymptotically as $\log(k)$, where $\gamma=0.5772156649\dots$ is the Euler-Mascheroni constant. The special values for even positive integer argument are given by the Euler's formula
\begin{equation}\label{eq:9}
\zeta(2k) = \frac{\mid B_{2k}\mid}{2(2k)!}(2\pi)^{2k},
\end{equation}
for which the value is expressed as a rational multiple of $\pi^{2k}$ where the constants $B_{2k}$ are Bernoulli numbers denoted such that $B_0=1$, $B_1=-1/2$, $B_2=1/6$ and so on. For odd positive integer argument, the values of this series converge to unique constants, which are not known to be expressed as a rational multiple of $\pi^{2k+1}$ as occurs in the even positive integer case.  For $n=3$, the value is commonly known as  Ap\'ery's constant, who proved its irrationality in 1979. Moreover, the Riemann zeta function has a prime product representation
\begin{equation}\label{eq:1}
\zeta(s) = \prod_{n=1}^\infty\ \left(1-\frac{1}{p^s_n}\right)^{-1},
\end{equation}
where $p_n$ is a sequence of $n$th prime numbers denoted such that $p_1=2$, $p_2 = 3$, $p_3=5$ and so on. In the previous article [4], we obtained an expression for the complex magnitude of Euler prime product as

\begin{equation}\label{eq:7}
\mid \zeta(\sigma+it) \mid^2 = \frac{\zeta(4\sigma)}{\zeta(2\sigma)}\prod_{n=1}^\infty \left(1-\frac{\cos(t\log p_n)}{\cosh(\sigma\log p_n)}\right)^{-1},
\end{equation}
for $\sigma>1$, which for a positive integer argument $\sigma=k$ simplifies the zeta terms using (3), resulting in
\begin{equation}\label{eq:11}
\mid \zeta(k+it) \mid=(2\pi)^{k}\sqrt{\frac{|B_{4k}|(2k)!}{|B_{2k}|(4k)!}}\prod_{n=1}^\infty \left(1-\frac{\cos(t\log p_n)}{\cosh(k\log p_n)}\right)^{-1/2}.
\end{equation}
Using this form, the first few special values of this representation are

\begin{equation}\label{eq:29}
\begin{aligned}
\zeta(1) &\sim \frac{\pi}{\sqrt{15}}\prod_{n=1}^k \left(1-\frac{2}{p_n+p_n^{-1}}\right)^{-1/2} \sim e^{\gamma}\log(p_k)\\
\zeta(2) &= \frac{\pi^2}{\sqrt{105}}\prod_{n=1}^\infty \left(1-\frac{2}{p_n^2+p_n^{-2}}\right)^{-1/2} \\
\zeta(3) &= \frac{\pi^3}{15} \sqrt{\frac{691}{3003}}\prod_{n=1}^\infty \left(1-\frac{2}{p_n^{3}+p_n^{-3}}\right)^{-1/2} \\
\zeta(4) &= \frac{\pi^4}{45} \sqrt{\frac{3617}{17017}}\prod_{n=1}^\infty \left(1-\frac{2}{p_n^{4}+p_n^{-4}}\right)^{-1/2} \\
\zeta(5) &= \frac{\pi^5}{225} \sqrt{\frac{174611}{323323}}\prod_{n=1}^\infty \left(1-\frac{2}{p_n^{5}+p_n^{-5}}\right)^{-1/2},
\end{aligned}
\end{equation}
where we let $t=0$ and reduced the hyperbolic cosine term. The value for $\zeta(1)$ in terms of Euler prime product representation is asymptotic to $e^{\gamma}\log(p_k)$ due to Merten's theorem [3,7]. The values for $\zeta(2k)$ have prime products of the form

\begin{equation}\label{eq:11}
\prod_{n=1}^\infty \left(1-\frac{2}{p_n^{2k}+p_n^{-2k}}\right)^{-1}\to \frac{a}{b}
\end{equation}
which yields a rational constant for some $a,b\in \mathbb{N}$ such that the full product (6) will simplify to Euler's formula (3). For example, for $k=2$ we have

\begin{equation}\label{eq:11}
\prod_{n=1}^\infty \left(1-\frac{2}{p_n^{2}+p_n^{-2}}\right)^{-1} = \frac{105}{36}.
\end{equation}
But for the odd case, however, it is not known whether the prime product simplifies to a rational constant. For example, for $k=3$ we have

\begin{equation}\label{eq:11}
\prod_{n=1}^\infty \left(1-\frac{2}{p_n^{3}+p_n^{-3}}\right)^{-1} = 1.46963883790849\dots.
\end{equation}

Also, the Riemann zeta function has an alternating series representation
\begin{equation}\label{eq:20}
\zeta(s) = \frac{1}{1-2^{1-s}}\sum_{n=1}^{\infty} \frac{(-1)^{n+1}}{n^s},
\end{equation}
whose domain of convergence is valid for $\Re(s)>0$, with some exceptions at $\Re(s)=1$ due to the constant factor. The non-trivial zeros of $\zeta(s)$ are constrained to a critical strip region $0<\sigma<1$, and hence they are zeros of the alternating series representation (11). Hardy proved that there is an infinity of non-trivial zeros on the critical line at $\sigma=1/2$, and the Riemann Hypothesis proposes that all non-trivial zeros should lie on the critical line at $\sigma=1/2$, which to date remains unproven [2]. For a positive integer index $q$, we denote a $q$th non-trivial zero on the critical line as $\rho'_q=1/2+it'_q$. The first few non-trivial zeros on the critical line have imaginary components  $t^{'}_1 = 14.13472514...$, $t^{'}_2 = 21.02203964...$, $t^{'}_3 = 25.01085758...$ (shown to eight decimal places) which are found numerically. From [5], we have an expression for the complex magnitude of this representation as
\begin{equation}\label{eq:20}
\mid \zeta(s) \mid^2 = C^2 \left[\sum_{n=1}^{\infty}\frac{1}{n^{2\sigma}}+2\sum_{n=1}^{\infty}\sum_{m=n+1}^{\infty}\frac{(-1)^{m}(-1)^{n}}{m^\sigma n^\sigma}\cos(t \log(m/n))\right],
\end{equation}
where $C$ is a constant and it is assumed that the index variables satisfy $m>n$ starting with $n=1$. If we then consider solutions to
\begin{equation}\label{eq:20}
\mid \zeta(s) \mid^2 = 0
\end{equation}
on the critical line at $\sigma=1/2$, then the complex magnitude induces an asymptotic harmonic series (2) when one considers the infinite series asymptotically as $k\to\infty$, from which we obtain the Euler-Mascheroni constant as
\begin{equation}\label{eq:20}
\gamma = \lim_{k\to \infty}\left(2\sum_{n=1}^{k}\sum_{m=n+1}^{k}\frac{(-1)^{m}(-1)^{n+1}}{\sqrt{mn}}\cos(t^{'}_q \log(m/n))-\log(k)\right)
\end{equation}
for all $q=1,2,3\dots$, which is expressed as a function of individual non-trivial zeros on the critical line.

On the other hand, a Dirichlet beta function is defined as
\begin{equation}\label{eq:1}
\beta(s) = \sum_{n=0}^{\infty}\frac{(-1)^n}{(2n+1)^s},
\end{equation}
which is convergent for $\Re(s)>0$, and for which the values at odd positive integer argument is expressed as a rational multiple of $\pi^{2k+1}$ as

\begin{equation}\label{eq:9}
\beta(2k+1) = \frac{|E_{2k}|}{2(2k)!}(\pi/2)^{2k+1},
\end{equation}
where the constants $E_{2k}$ are Euler numbers denoted such that $E_0=1$, $E_2=-1$, $E_4=5$ and so on. More specific details about Bernoulli and Euler numbers can be found in reference [1]. The first few special values of the Dirichlet beta function are

\begin{equation}\label{eq:9}
\begin{aligned}
\beta(1) &=\frac{\pi}{4} \\
\beta(2) &=0.915965594177219\dots \\
\beta(3) &=\frac{\pi^3}{32} \\
\beta(4) &=0.988944551741105\dots \\
\beta(5) &=\frac{5\pi^5}{1536},
\end{aligned}
\end{equation}
in which the even order case admits a unique value that is not known to be expressed as a rational multiple of $\pi^{2k}$. The value for $\beta(1)$ is related to $\tan^{-1}(1)$. The value for $\beta(2)$ is also known as the Catalan's constant $G$, but it is not known whether it is irrational. Also, the Dirichlet beta function has a prime product representation as

\begin{equation}\label{eq:1}
\beta(s) = \prod_{n=2}^\infty\ \left(1-\frac{\chi_4(p_n)}{p^s_n}\right)^{-1},
\end{equation}
which is convergent for $\Re(s)\geq 1$, where $\chi_4(n)$ is the Dirichlet character modulo $4$ defined as

\begin{equation}
\chi_4(n)= \left \{
\begin{aligned}
&1, &&\text{if}\ n\equiv 1 \text{ mod } 4 \\
-&1,&& \text{if}\ n\equiv 3 \text{ mod } 4 \\
&0, && \text{if}\ n \text{ is even},
\end{aligned} \right.
\end{equation}
but for a prime argument $p>2$, it reduces to a simple alternating sign function

\begin{equation}\label{eq:1}
\chi_4(p) = (-1)^{\frac{p-1}{2}}.
\end{equation}
In this article, we will explore the complex magnitude of the prime product representation of the Dirichlet beta function in more detail and develop formulas for $\beta(k)$ up to $11$th order, including the Catalan's constant, which is summarized in Appendix A. Also, the numerical computation of these results is summarized in Appendix B.

Furthermore, the Dirichlet beta function has non-trivial zeros in the critical strip. We do not know whether there is an infinity of zeros on the critical line, but it most likely is just as in the $\zeta(s)$ case. We also do not know whether all zeros lie on the critical line. The zeros $\beta(s)$ found numerically thus far do indeed lie on the critical line, so we denote a $q$th non-trivial zero as $\rho'_{Bq}=1/2+ir'_q$. The first few non-trivial zeros on the critical line have imaginary components  $r^{'}_1 = 6.02094890...$, $r^{'}_2 = 10.24377030...$, $r^{'}_3 = 12.98809801...$ (shown to eight decimal places) which are found numerically. Hence, we find an expression for the complex magnitude of the alternating series representation (15) and then similarly obtain a formula for $\gamma$ as a function of non-trivial zeros of the Dirichlet beta function on the critical line.

Finally, using the presented results, we will investigate the asymptotic behavior of the Euler prime product on the critical line, and note that it contains the $\zeta(1)$ singularity.

\section{Dirichlet beta prime product}
By substituting the complex argument to the beta function, and following the derivation in [4], we obtain a similar result

\begin{equation}\label{eq:3}
\beta(\sigma+it) = \prod_{n=1}^\infty\ \frac{p_n^{\sigma}-\chi_{4}(p_n)p_n^{it}}{p_n^{\sigma}+\chi_{4}^{2}(p_n)p_n^{-\sigma}-2\chi_{4}(p_n)\cos{(t\log{p_n}})},
\end{equation}
hence we obtain a formula for the complex magnitude as

\begin{equation}\label{eq:7}
\mid \beta(\sigma+it) \mid^2 = \prod_{n=1}^\infty \left(1+\frac{\chi_{4}(p_n)^2}{p_n^{2\sigma}}\right)^{-1}\left(1-\chi_4(p_n)\frac{\cos(t\log p_n)}{\cosh(\sigma\log p_n)}\right)^{-1}.
\end{equation}
Since $\chi_4(2)=0$ and $\chi_4(p)^2=1$ for $p>2$, we can extract the term $(1+1/2^{2\sigma})$ to obtain the final representation as
\begin{equation}\label{eq:7}
\mid \beta(\sigma+it) \mid^2 = \left(1+\frac{1}{2^{2\sigma}}\right)\frac{\zeta(4\sigma)}{\zeta(2\sigma)}\prod_{n=2}^\infty \left(1-\chi_4(p_n)\frac{\cos(t\log p_n)}{\cosh(\sigma\log p_n)}\right)^{-1}
\end{equation}
for $\sigma \geq 1$. By simplifying the zeta terms further, and letting $\sigma$ be a positive integer $k$, then we have

\begin{equation}\label{eq:7}
\mid \beta(k+it) \mid = \pi^{k}\sqrt{\left(1+2^{2k}\right)\frac{|B_{4k}|(2k)!}{|B_{2k}|(4k)!}}\prod_{n=2}^\infty \left(1-\chi_4(p_n)\frac{\cos(t\log p_n)}{\cosh(k\log p_n)}\right)^{-1/2}.
\end{equation}
This form is a slight modification of (6), having only a $1+2^{2k}$ in the factor and the Dirichlet character in the prime product. From this, we can observe that the prime product term at odd positive integer argument yields a rational constant of the form
\begin{equation}\label{eq:7}
\prod_{n=2}^\infty \left(1-\chi_4(p_n)\frac{2}{p_n^{2k+1}+p_n^{-(2k+1)}}\right)^{-1}=\frac{a}{b}
\end{equation}
for some $a,b\in \mathbb{N}$ such that the full product (24) will simplify to (16). We next will evaluate this representation for the first few special values.

\section{Evaluation of $\beta(1)$}
By setting $k=1$ and $t=0$ into equation (24) yields

\begin{equation}\label{eq:7}
\beta(1) = \frac{\pi}{2\sqrt{3}}\prod_{n=2}^\infty \left(1-\chi_4(p_n)\frac{2}{p_n+p_n^{-1}}\right)^{-1/2}.
\end{equation}
From this, we can deduce that the prime product term yields a rational constant
\begin{equation}\label{eq:7}
\prod_{n=2}^\infty \left(1-\chi_4(p_n)\frac{2}{p_n+p_n^{-1}}\right)^{-1} = \frac{3}{4},
\end{equation}
where we obtain the familiar result

\begin{equation}\label{eq:7}
\beta(1) = \frac{\pi}{4}.
\end{equation}
In the next example, we evaluate the complex magnitude of $\beta(1+i)$ as

\begin{equation}\label{eq:7}
\mid \beta(1+i) \mid = \frac{\pi}{2\sqrt{3}}\prod_{n=2}^\infty \left(1-\chi_4(p_n)\frac{\cos(\log p_n)}{\cosh(\log p_n)}\right)^{-1/2} = 0.879512157843\dots.
\end{equation}

\section{Evaluation of $\beta(2)$}
By setting $k=2$ and $t=0$ into equation (24) yields a formula for Catalan's constant

\begin{equation}\label{eq:7}
\beta(2) = \frac{\pi^2}{4} \sqrt{\frac{17}{105}}\prod_{n=2}^\infty \left(1-\chi_4(p_n)\frac{2}{p_n^2+p_n^{-2}}\right)^{-1/2}.
\end{equation}
The prime product term results in a constant
\begin{equation}\label{eq:7}
\prod_{n=2}^\infty \left(1-\chi_4(p_n)\frac{2}{p_n^2+p_n^{-2}}\right)^{-1} = 0.85117565043382,
\end{equation}
which is not known to be rational, as in the odd order case. In the next example, we evaluate the complex magnitude of $\beta(2+i)$ as

\begin{equation}\label{eq:7}
\mid \beta(2+i) \mid =\frac{\pi^2}{4}\sqrt{\frac{17}{105}}\prod_{n=2}^\infty \left(1-\chi_4(p_n)\frac{\cos(\log p_n)}{\cosh(2\log p_n)}\right)^{-1/2} = 0.954186328286\dots.
\end{equation}

\section{Evaluation of $\beta(3)$}
By setting $k=3$ and $t=0$ into equation (24) yields

\begin{equation}\label{eq:7}
\beta(3) = \frac{\pi^3}{24}\sqrt{\frac{691}{1155}}\prod_{n=2}^\infty \left(1-\chi_4(p_n)\frac{2}{p_n^3+p_n^{-3}}\right)^{-1/2}.
\end{equation}
From this, we can deduce that the prime product term yields a rational constant
\begin{equation}\label{eq:7}
\prod_{n=2}^\infty \left(1-\chi_4(p_n)\frac{2}{p_n^{3}+p_n^{-3}}\right)^{-1} = \frac{10395}{11056},
\end{equation}
where we obtain the familiar result

\begin{equation}\label{eq:7}
\beta(3) = \frac{\pi^3}{32}.
\end{equation}
In the next example, we evaluate the complex magnitude of $\beta(3+i)$ as

\begin{equation}\label{eq:7}
\mid \beta(3+i) \mid =\frac{\pi^3}{24} \sqrt{\frac{691}{1155}}\prod_{n=2}^\infty \left(1-\chi_4(p_n)\frac{\cos(\log p_n)}{\cosh(3\log p_n)}\right)^{-1/2} = 0.983806269077\dots.
\end{equation}

These formulas are summarized up to $11$th order in Appendix A and are numerically computed in Appendix B.

\section{Dirichlet beta alternating series}
We consider evaluating the complex magnitude of an alternating series representation of the Dirichlet beta function (15) for $s=\sigma+it$ on the critical strip $0<\sigma<1$ which results in a form

\begin{equation}\label{eq:20}
\mid \beta(s) \mid^2 = A^2+B^2,
\end{equation}
where the constants $A$ and $B$ are the real and imaginary parts of the infinite sum term of (15) as

\begin{equation}\label{eq:20}
\begin{aligned}
A =  \sum_{n=0}^{\infty} \frac{(-1)^{n}}{(2n+1)^\sigma}\cos(t\log(2n+1)) \\
B =  -\sum_{n=0}^{\infty} \frac{(-1)^{n}}{(2n+1)^\sigma}\sin(t\log(2n+1)).
\end{aligned}
\end{equation}
We note that since $A$ and $B$ are convergent in the critical strip, their squares are also convergent, and so the complex magnitude is convergent. When one expands the sum of squares of $A$ and $B$, we obtain a double series formed by the complex magnitude, and using a trigonometric identity simplifies the sum of cosine and sine products into a more compact form

\begin{equation}\label{eq:20}
\mid \beta(s) \mid^2 = \sum_{m=0}^{\infty}\sum_{n=0}^{\infty} \frac{(-1)^{m}(-1)^{n}}{(2m+1)^\sigma (2n+1)^\sigma}\cos(t (\log(2m+1)-\log(2n+1))).
\end{equation}
We then observe that due to the symmetry of the double series we can re-arrange the terms into diagonal ($m=n$) and off-diagonal ($m\neq n$) sums, hence we have

\begin{equation}\label{eq:20}
\mid \beta(s) \mid^2 = \sum_{n=0}^{\infty}\frac{1}{(2n+1)^{2\sigma}}+2\sum_{n=0}^{\infty}\sum_{m=n+1}^{\infty}\frac{(-1)^{m}(-1)^{n}}{(2m+1)^\sigma (2n+1)^\sigma}\cos\left(t \log\frac{2m+1}{2n+1}\right),
\end{equation}
where it is assumed the index variables $m$ and $n$ are positive integers starting with $n=0$ and satisfying $m>n$. Now we consider solutions to

\begin{equation}\label{eq:20}
\mid \beta(s) \mid^2 = 0
\end{equation}
on the critical strip, which are the zeros of the Dirichlet beta function, since the real and imaginary parts are coupled together by the complex magnitude. Hence, it is evident from (40) that if $\sigma=1/2$ implies the solutions must satisfy an equation

\begin{equation}\label{eq:20}
\lim_{k\to\infty}\left[\sum_{n=0}^{k}\frac{1}{(2n+1)}+2\sum_{n=0}^{k}\sum_{m=n+1}^{k}\frac{(-1)^{m}(-1)^{n}}{\sqrt{(2m+1)(2n+1)}}\cos\left(r'
_q \log\frac{2m+1}{2n+1}\right)\right]=0,
\end{equation}
where $r^{'}_q$ is the $q$th imaginary component of a non-trivial zero of $\beta(s)$ on the critical line. We note that these sums are divergent, but they must cancel since (39) is convergent. Therefore, when we treat these divergent sums in a limiting sense, then we can deduce an asymptotic relationship for harmonic series for an odd positive integer index as
\begin{equation}\label{eq:20}
\sum_{n=0}^{k}\frac{1}{(2n+1)} \sim \frac{1}{2}\gamma+\log(2)+\frac{1}{2}\log(k),
\end{equation}
which easily follows from (2). The second term is the new series expressed in terms of a single zero of $\beta(s)$ on the critical line. It then follows that the formula for the Euler-Mascheroni constant is

\begin{equation}\label{eq:20}
\gamma = \lim_{k\to \infty}\left[4\sum_{n=0}^{k}\sum_{m=n+1}^{k}\frac{(-1)^{m}(-1)^{n+1}}{\sqrt{(2m+1)(2n+1)}}\cos\left(r'_q \log\frac{2m+1}{2n+1}\right)-\log(4k)\right]
\end{equation}
as $k\to\infty$, where in this form, we absorbed the negative sign by the alternating sign index $n+1$. The index $q$ of a $r'_q$ zero is possibly valid from $1$ to infinity. It is not known whether there is an infinity of zeros of the Dirichlet beta function on the critical line to the best of my knowledge. In Table 1, we numerically compute the Euler-Mascheroni constant by equation (44) for the first ten zeros at $k=10^5$ in the Matlab software package. We used the non-trivial zeros from the LMFDB database [6].  We note that convergence is stable and accurate to four decimal places. The number of off-diagonal elements of the double sum is $k^2/2-k$, so for $k=10^5$, that is almost $5$ billion computations, which may take several hours to compute on a present-day desktop computer.

\begin{table}[ht]
\caption{Evaluation of $\gamma$ by Equation (44) for $k=10^5$}
\centering 
\begin{tabular}{c c c c} 
\hline\hline 
$q$ & $r_q^{'}$ & $\gamma$ Eq.(44)\\ [0.5ex] 
\hline 
$1$  & 6.0209489046976 & 0.577223164876768  \\
$2$  & 10.2437703041665 & 0.577223164850274 \\
$3$  & 12.9880980123124 & 0.577223164851997  \\
$4$  & 16.3426071045872 & 0.577223164872144  \\
$5$  & 18.2919931961235 & 0.577223164930105  \\
$6$  & 21.4506113439835 & 0.577223164868668 \\
$7$  & 23.2783765204595 & 0.577223164880770  \\
$8$  & 25.7287564250887 & 0.577223164815470 \\
$9$  & 28.3596343430253 & 0.577223164835386 \\
$10$ & 29.6563840145932 & 0.577223164855655
\\ [1ex] 
\hline 
\end{tabular}
\label{table:nonlin} 
\end{table}

\section{Euler prime product on the critical line}
The Euler prime product converges for $\Re(s)>1$ or $\Re(s)\geq 1$ in case of the Dirichlet beta function, and so it is divergent in the critical strip. Since $\chi^2=1$, we note that the internal structure of the Euler prime product (5) on the critical line contains the $\zeta(1)$ singularity. We can use Merten's theorem as an asymptotic representation which for the Riemann zeta function leads to

\begin{equation}\label{eq:7}
\mid\zeta_{P}(1/2+it) \mid^2 \sim \frac{\pi^2}{6}\frac{1}{e^{\gamma}\log(p_k)}\prod_{n=1}^k \left(1-\frac{\cos(t\log p_n)}{\cosh(1/2\log p_n)}\right)^{-1}
\end{equation}
where we distinguish $\zeta(s)$ with $\zeta_P(s)$ to indicate a prime product representation in the critical strip, which is not an analytical continuation of $\zeta(s)$. Similarly for the Dirichlet beta function we have
\begin{equation}\label{eq:7}
\mid \beta_P(1/2+it) \mid^2 \sim \frac{\pi^2}{4}\frac{1}{e^{\gamma}\log(p_k)}\prod_{n=1}^k \left(1-\chi_4(p_n)\frac{\cos(t\log p_n)}{\cosh(1/2\log p_n)}\right)^{-1}.
\end{equation}
We note that the first term on the left is asymptotic to $\log(p_k)$ due to the $\zeta(1)$ singularity, however, the second term involving products of $\cos(t\log(p_n))$ is highly oscillatory. In Figure 1 we plot $\mid\zeta_P(1/2+it)\mid^2$ and in Figure 2 we plot $\mid\beta_P(1/2+it)\mid^2$ for $k=10^3$ and compare with the analytically continued representation. The plots become highly oscillatory and erratic, especially when $k$ is increased, but also we observe that nulls form at locations of non-trivial zeros. When one essentially considers the sum of cosine waves with $\log(n)$ angular frequencies such as by equation (12), then these nulls form as localized nodes of standing waves. In the case of the Euler prime product representation, these nulls may not approach zero because of the divergence of the second term. In case of an analytically continued representation for $\zeta(s)$ or $\beta(s)$ to the critical strip, such as by equation (11) and (15) respectively, then these nulls approach zero as the non-trivial zeros.

Finally, instead of using the $\log(p_n)$ from Merten's theorem as an asymptotic representation of $\zeta(1)$ and if instead, we use the asymptotic representation in equation (7), then we can obtain the complex magnitude for equations (45) and (46) in terms of $\sqrt{\pi}$ as

\begin{equation}\label{eq:7}
\mid\zeta_{P}(1/2+it) \mid \sim \sqrt{\pi}\left(\frac{5}{12}\right)^{1/4}\prod_{n=1}^k \left[\frac{\sqrt{p_n^2+1}}{p_n-1}\left(1-\frac{\cos(t\log p_n)}{\cosh(1/2\log p_n)}\right)\right]^{-1/2}
\end{equation}
and similarly, for the Dirichlet beta function, we have

\begin{equation}\label{eq:7}
\mid\beta_{P}(1/2+it) \mid \sim \frac{\sqrt{\pi}}{2}\left(15\right)^{1/4}\prod_{n=1}^k \left[\frac{\sqrt{p_n^2+1}}{p_n-1}\left(1-\chi_4(p_n)\frac{\cos(t\log p_n)}{\cosh(1/2\log p_n)}\right)\right]^{-1/2}.
\end{equation}

\begin{figure}[h]
  \centering
  \includegraphics[width=100mm]{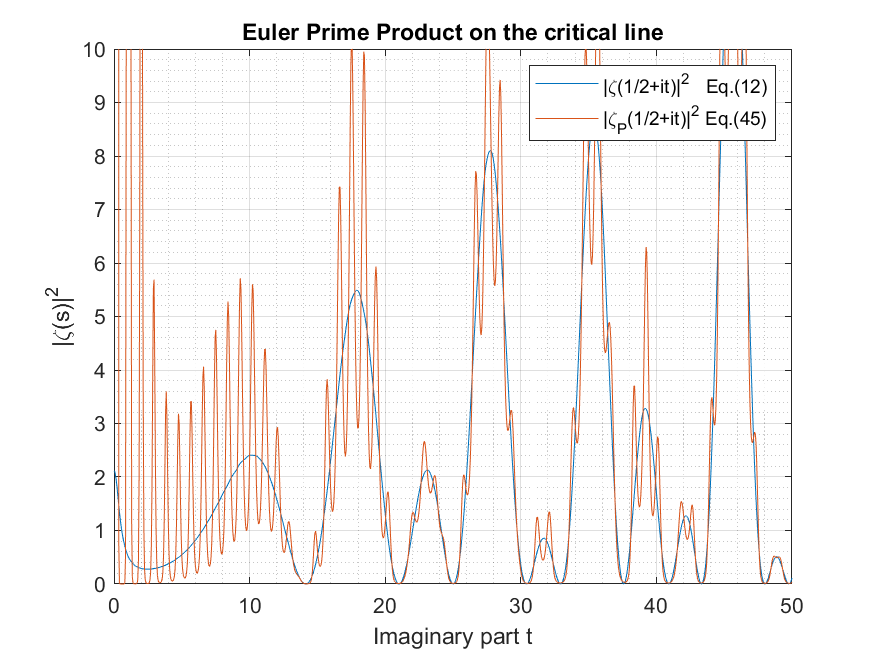}\\
  \caption{Euler prime product on the critical line for $|\zeta(s)|^2$ for $k=10^3$.}\label{1}
\end{figure}

\begin{figure}[h]
  \centering
  \includegraphics[width=100mm]{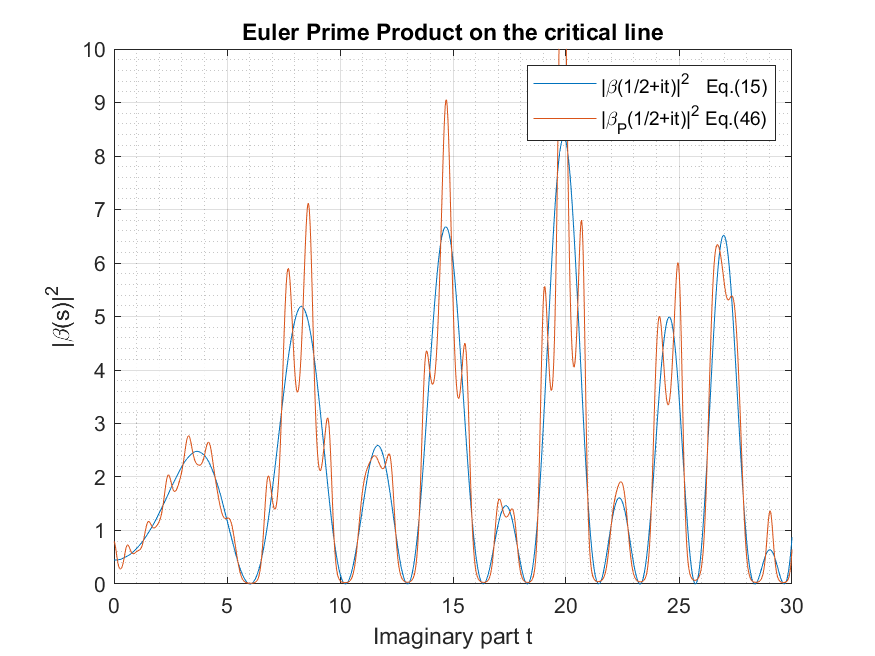}\\
  \caption{Euler prime product on the critical line for $|\beta(s)|^2$ for $k=10^3$.}\label{1}
\end{figure}

\section{Conclusion}
We obtained an expression for the complex magnitude of the Euler prime product representation of the Dirichlet beta function, from which we can generate identities involving $\pi^k$ term for an even and odd positive integer argument. We also obtained an expression for the complex magnitude of the alternating series representation (15), from which we can deduce an expression for the Euler-Mascheroni constant as a function of non-trivial zeros of the Dirichlet beta function on the critical line, which is induced by the solutions of $|\beta(s)|^2=0$. We also numerically validate the presented results and summarize them in the Appendix section. We also investigated the Euler prime product on the critical line and noted that it contains the $\zeta(1)$ singularity. It also contains prime product terms that are highly oscillatory due to a sum of cosine waves with $\log(p_n)$ angular frequencies, which at the location of non-trivial zeros appear as nodes of standing waves. The Euler prime product doesn't vanish exactly at the non-trivial zeros, but the null is very strong (at least when one approximates a partial Euler product up to $k$) where the amplitudes of the cosine waves in the oscillatory product term diverge. Finally, we note that in all of the investigated cases, the complex magnitude induces the $\zeta(1)$ singularity on the critical line, which contains the Euler-Mascheroni constant due to either the harmonic series or in case of the Euler prime product, the Merten's theorem.

\texttt{Email: art.kawalec@gmail.com}

\clearpage
\section{Appendix A}
The evaluation of equation (24) in the Mathematica software package up to the $11$th order:

\begin{equation}\label{eq:29}
\beta(1) = \frac{\pi}{2\sqrt{3}}\prod_{n=2}^\infty \left(1-\chi_4(p_n)\frac{2}{p_n+p_n^{-1}}\right)^{-1/2}
\nonumber
\end{equation}

\begin{equation}\label{eq:29}
\beta(2) = \frac{\pi^2}{4} \sqrt{\frac{17}{105}}\prod_{n=2}^\infty \left(1-\chi_4(p_n)\frac{2}{p_n^{2}+p_n^{-2}}\right)^{-1/2}
\nonumber
\end{equation}

\begin{equation}\label{eq:30}
\beta(3) = \frac{\pi^3}{24} \sqrt{\frac{691}{1155}}\prod_{n=2}^\infty \left(1-\chi_4(p_n)\frac{2}{p_n^{3}+p_n^{-3}}\right)^{-1/2}
\nonumber
\end{equation}

\begin{equation}\label{eq:31}
\beta(4) = \frac{\pi^4}{720} \sqrt{\frac{929569}{17017}}\prod_{n=2}^\infty \left(1-\chi_4(p_n)\frac{2}{p_n^{4}+p_n^{-4}}\right)^{-1/2}
\nonumber
\end{equation}

\begin{equation}\label{eq:32}
\beta(5) = \frac{\pi^5}{1440} \sqrt{\frac{7159051}{323323}}\prod_{n=2}^\infty \left(1-\chi_4(p_n)\frac{2}{p_n^{5}+p_n^{-5}}\right)^{-1/2}
\nonumber
\end{equation}

\begin{equation}\label{eq:33}
\beta(6) = \frac{\pi^6}{20160} \sqrt{\frac{56963745931}{129543843}}\prod_{n=2}^\infty \left(1-\chi_4(p_n)\frac{2}{p_n^{6}+p_n^{-6}}\right)^{-1/2}
\nonumber
\end{equation}

\begin{equation}\label{eq:34}
\beta(7) = \frac{\pi^7}{1814400} \sqrt{\frac{383384156611}{1062347}}\prod_{n=2}^\infty \left(1-\chi_4(p_n)\frac{2}{p_n^{7}+p_n^{-7}}\right)^{-1/2}
\nonumber
\end{equation}

\begin{equation}\label{eq:35}
\beta(8) = \frac{\pi^8}{3628800} \sqrt{\frac{505245773078238529}{3454415680001}}\prod_{n=2}^\infty \left(1-\chi_4(p_n)\frac{2}{p_n^{8}+p_n^{-8}}\right)^{-1/2}
\nonumber
\end{equation}

\begin{equation}\label{eq:36}
\beta(9) = \frac{\pi^9}{79833600} \sqrt{\frac{2868364599282829033657}{399908500933305}}\prod_{n=2}^\infty \left(1-\chi_4(p_n)\frac{2}{p_n^{9}+p_n^{-9}}\right)^{-1/2}
\nonumber
\end{equation}

\begin{equation}\label{eq:37}
\beta(10) = \frac{\pi^{10}}{6227020800} \sqrt{\frac{16103843159579478297227731}{3642221497651835}}\prod_{n=2}^\infty \left(1-\chi_4(p_n)\frac{2}{p_n^{10}+p_n^{-10}}\right)^{-1/2}
\nonumber
\end{equation}

\begin{equation}\label{eq:38}
\beta(11) = \frac{\pi^{11}}{12454041600} \sqrt{\frac{2122567668414730590074148073}{1184508331617858905}}\prod_{n=2}^\infty \left(1-\chi_4(p_n)\frac{2}{p_n^{11}+p_n^{-11}}\right)^{-1/2}.
\nonumber
\end{equation}

\newpage
\section{Appendix B}

In the following table, we summarize the numerical computation of equation (24) in the Mathematica software package to $15$ decimal places for the first $1000$ prime product terms. We also note that convergence is slower for smaller arguments and that all the higher-order argument converged faster.

\begin{table}[ht]
\caption{Evaluation of $\beta(k)$ and the new formula for $\beta(k)$ respectively} 
\centering 
\begin{tabular}{c c c} 
\hline\hline 
k & $\beta(k)$ & $\beta(k)$ Equation (24) \\ [0.5ex] 
\hline 
1  & 0.785398163397448 & 0.786450346758764 \\ 
2  & 0.915965594177219 & 0.915965663664403\\
3  & 0.968946146259369  & 0.968946146264864 \\
4  & 0.988944551741105   & 0.988944551741106 \\
5  & 0.996157828077088 & 0.996157828077088 \\
6  & 0.998685222218438 & 0.998685222218438 \\
7  & 0.999554507890540 & 0.999554507890540 \\
8  & 0.999849990246830  & 0.999849990246830 \\
9  &  0.999949684187220  & 0.999949684187220 \\
10  & 0.999983164026197 & 0.999983164026197 \\
11 & 0.999994374973824  & 0.999994374973824
 \\ [1ex] 
\hline 
\end{tabular}
\label{table:nonlin} 
\end{table}

\end{document}